\documentclass[10pt]{article}%
\usepackage{amsmath}
\usepackage{graphicx}
\usepackage{amsfonts}
\usepackage{amssymb}
\usepackage{times}%

\begin{document}

\begin{center}
.

\vspace{0.4in}

{\LARGE \textbf{Discrete Dynamical Modeling and Analysis of the \textit{R-S}
Flip-Flop Circuit}}

\vspace{0.3in}

\smallskip\textsf{Denis Blackmore}

\textsf{Department of Mathematical Sciences and}

\textsf{Center for Applied Mathematics and Statistics}

\textsf{New Jersey Institute of Technology}

\textsf{Newark, NJ 07102-1982}

\textsf{deblac@m.njit.edu}

$\ast\;\ast\;\ast$

\smallskip\textsf{Aminur Rahman}

\textsf{Department of Mathematical Sciences}

\textsf{New Jersey Institute of Technology}

\textsf{Newark, NJ 07102-1982}

\textsf{ar276@njit.edu}

$\ast\;\ast\;\ast$\textsf{\ }

\textsf{Jigar Shah }

\textsf{Department of Computer and Electrical Engineering}

\textsf{New Jersey Institute of Technology}

\textsf{Newark, NJ 07102-1982}

\textsf{jds22@njit.edu }
\end{center}

\vspace{0.2in}

\noindent\textbf{ABSTRACT: }A simple discrete planar dynamical model for the
ideal (logical) \textit{R-S} flip-flop circuit is developed with an eye toward
mimicking the dynamical behavior observed for actual physical realizations of
this circuit. It is shown that the model exhibits most of the qualitative
features ascribed to the \textit{R-S} flip-flop circuit, such as an intrinsic
instability associated with unit set and reset inputs, manifested in a chaotic
sequence of output states that tend to oscillate among all possible output
states, and the existence of periodic orbits of arbitrarily high period that
depend on the various intrinsic system parameters. The investigation involves
a combination of analytical methods from the modern theory of discrete
dynamical systems, and numerical simulations that illustrate the dazzling
array of dynamics that can be generated by the model. Validation of the
discrete model is accomplished by comparison with certain Poincar\'{e} map
like representations of the dynamics corresponding to three-dimensional
differential equation models of electrical circuits that produce
\textit{R}-\textit{S} flip-flop behavior.

\bigskip


\medskip


\bigskip

\section{Introduction}

The ideal \emph{R-S flip-flop circuit} is a logical feedback circuit that can
be described most efficiently in terms of Fig. 1, with input/output behavior
described in Table 1, which shows the \emph{set} ($S$) and \emph{reset} ($R$)
inputs feeding into the simple circuit comprised of two \emph{nor gates} and
the corresponding outputs $Q$ and $P$ that are generated. The input to this
circuit may be represented as $(S,R)$ and the output by $(Q,P)$, so the
association of the input to the output, denoted by $(S,R)\rightarrow(Q,P)$ may
be regarded as the action of a map from the plane $\mathbb{R}^{2}%
:=\{(x,y):x,y\in\mathbb{R}\}$ into itself, where $\mathbb{R}$ denotes the real
numbers. Our goal, from this mapping perspective, is to construct a simple
nonlinear map of the plane that models the logical properties of the
\textit{R-S} flip-flop circuit, with iterates (discrete dynamics) that at
least qualitatively capture most of its interesting dynamics, both those that
are intuitive and those that have been observed in studies of physical circuit
simulations - especially in certain critical cases that we shall describe in
the sequel.

\begin{figure}[th]
\centering
\includegraphics[width=3.5in]{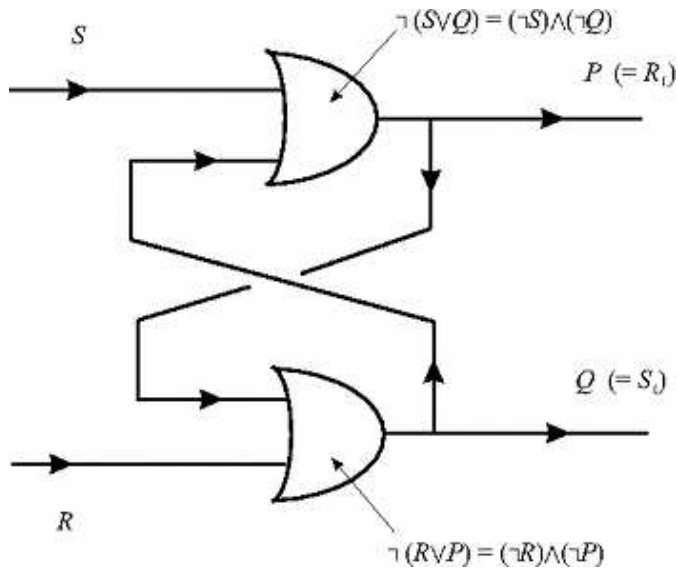}\caption{$R$-$S$
flip-flop circuit}%
\label{circuit}%
\end{figure}

%\begin{figure}[th]
%\centering
%\includegraphics[width=6in]{vbchaos.jpg}
%\caption{R-S flip-flop circuit}
%\label{circuit1}
%\end{figure}

\noindent The binary input/output behavior, with 0 and 1 representing false
and true, respectively, is given in the following table.

\begin{center}%
\begin{tabular}
[c]{|c|c|c|c|}\hline
$S$ & $R$ & $S_{1}:=Q$ & $R_{1}:=P$\\\hline
1 & 0 & \multicolumn{1}{|r|}{$%
\begin{array}
[c]{c}%
1
\end{array}
$} & \multicolumn{1}{|r|}{$%
\begin{array}
[c]{c}%
0
\end{array}
$}\\\hline
0 & 1 & \multicolumn{1}{|r|}{$%
\begin{array}
[c]{c}%
0
\end{array}
$} & \multicolumn{1}{|r|}{$%
\begin{array}
[c]{c}%
1
\end{array}
$}\\\hline
1 & 1 & \multicolumn{1}{|r|}{$%
\begin{array}
[c]{c}%
0
\end{array}
$} & \multicolumn{1}{|r|}{$%
\begin{array}
[c]{c}%
0
\end{array}
$}\\\hline
0 & 0 & \multicolumn{1}{|r|}{or $\
\begin{array}
[c]{c}%
1\\
0
\end{array}
$} & \multicolumn{1}{|r|}{or $%
\begin{array}
[c]{c}%
0\\
1
\end{array}
$}\\\hline
\end{tabular}

\medskip

Table 1. Binary input/output of $R$-$S$ flip-flop circuit

\medskip
\end{center}

Note that the input/output behavior is not well defined when the set and reset
values are both $0$; a situation that can be remedied by identifying $(1,0)$
with $(0,1)$. Such an identification strategy is consistent with the obvious
symmetry of the circuit with respect to $S$ and $R$, and in fact we shall
employ this in the sequel when we describe and analyze our $R$-$S$ flip-flop
map model. If we make this symmetry identification, it leads to an unambiguous
input/output behavior for the circuit, and from an abstract perspective, it is
not immediately clear why the state $(1,1)$ should be problematical - other
than that it is the only state that produces non-complementary output values
$0,0$. Notwithstanding the well defined behavior of the abstract
\textit{R}-\textit{S} flip-flop circuit when the symmetry identification is
imposed, actual physical circuit models comprised of elements such as
capacitors, diodes, inductors and resistors exhibit highly oscillatory, very
unstable and even chaotic dynamics (\emph{metastable operation}), as
experimentally observed in such studies as \cite{Chaney, KA, kac, LMH}. This
type of highly irregular dynamical behavior has also been found in
\textit{R}-\textit{S} flip-flop realizations in the context of Josephson
junctions \cite{ztkbn}. In general, it is found that physical models of
\textit{R}-\textit{S} flip-flop circuits invariably generate very complex
dynamics belying the simplicity of the abstract logical circuit, which can
plausibly be ascribed to the fact that every real model of this circuit must
have inherent time sequencing characteristics due to the finite speed of
electromagnetic waves.

Continuous dynamical systems representations of \textit{R}-\textit{S}
flip-flop circuits derived from the usual circuit equations applied to the
physical models, typically lead to three-dimensional systems of first order,
autonomous, (possibly discontinuous) piecewise smooth (and usually piecewise
linear) differential equations such as those investigated by Murali \emph{et
al}. \cite{msd}, Okazaki \emph{et al}. \cite{okt} and Ruzbehani \emph{et al}.
\cite{rzw}. These mathematical models are traceable back to the pioneering
work of Moser \cite{Moser}, and are subsumed by the famous circuit of Chua and
its generalizations \cite{chua, CWHZ}. Experimental observations on real
circuit models of \textit{R}-\textit{S} flip-flops, together with numerical
simulations utilizing such tools as SPICE \cite{hamill}, and analytical
studies employing the tools of modern dynamical systems theory (see
\emph{e.g}. \cite{guho, kathas, palmel, wigbook}) including Poincar\'{e}
sections and Melnikov functions such as in \cite{CWHZ, Danca, kac, okt, rzw}
have painted a rather compelling picture of the extreme complexity of the
dynamical possibilities.

Naturally, when one has a series of quite successful mathematical
representations of phenomena or processes as is the case for \textit{R-S}
flip-flop circuit behavior, it leads to the question of the possible existence
of simpler models. One cannot expect to reduce the dimension of the continuous
dynamical models, since it is impossible for two-dimensional systems of
autonomous differential equations to have chaotic solutions. However,
two-dimensional discrete dynamical systems are an enticing possibility since
they can exhibit almost all types of complex dynamics - including chaotic
regimes. Moreover, there have been some studies, albeit just a few, including
those of Danca \cite{Danca} and Hamill \emph{et al}. \cite{hdj}, which give
strong indications of the potential of modeling circuits such as the
\textit{R}-\textit{S} flip-flop using two-dimensional difference equations or
iterated planar maps. Encouraged by this literature on discrete dynamical
models and relying heavily on our knowledge of dynamical systems theory and
physical intuition, we have developed a discrete - essentially
phenomenological -dynamical model generated by iterates of a rather simple
nonlinear (quadratic), two-parameter planar map, which we present and analyze
in this paper.

Our investigation is organized as follows. First, in Section 2, we define our
simple planar map model - with iterates producing the dynamic behavior that we
shall show mimics that observed in real \textit{R}-\textit{S} flip-flop
circuits. Moreover, we derive some basic properties of the map related to
fixed points and the existence of a local inverse. This is followed in Section
3 with a more thorough analysis of the fixed points of the map - including a
local stability analysis and an analysis of stable and unstable manifolds -
under very mild, and quite reasonable, restrictions on the two map parameters.
Next, in Section 4, we fix the value of one of the parameters and prove the
existence of a Hopf bifurcation at an interior fixed point as the other
parameter is varied. We also note what appears to be a kind of Hopf
bifurcation cascade (manifested by an infinite sequence) as the parameter is
increased, which suggests the existence of extreme oscillatory behavior
culminating in chaos. We then prove the existence of chaotic dynamics in
Section 5. This is followed in Section 6 by a comparison of the dynamics of
our model with other results in the literature, primarily from the perspective
of planar Poincar\'{e} maps. Naturally, liberal use is made in this section of
numerical simulations of our model for comparison purposes. Finally, in
Section 7, we summarize and underscore some of the more important results of
this investigation, and identify several interesting directions for future
related research.

\section{A Discrete Model}

Based upon the ideal mathematical properties of the \textit{R-S} flip-flop
circuit, which were briefly delineated in the preceding section, and a
knowledge of the interesting dynamical characteristics of actual physical
circuits constructed to perform like the \textit{R-S} flip-flop circuit (see
\emph{e.g}. \cite{hdj, kac, okt, rzw, ztkbn}, and also \cite{Danca} for
related results) , we postulate the following simple, quadratic, two-parameter
planar map model:%

\begin{equation}
\Phi=\left(  \varphi,\psi\right)  =\Phi_{\lambda,\mu}=\left(  \varphi
_{\lambda,\mu},\psi_{\lambda,\mu}\right)  :\mathbb{R}^{2}\rightarrow
\mathbb{R}^{2}, \label{e1}%
\end{equation}
where $\lambda$ and $\mu$ are positive parameters, and the coordinate
functions $\varphi$ and $\psi$ are defined as%
\begin{align}
\varphi\left(  x,y\right)   &  =\varphi_{\lambda,\mu}\left(  x,y\right)
:=1-x\left[  \lambda\left(  1-x\right)  +y\right]  ,\nonumber\\
\psi\left(  x,y\right)   &  =\psi_{\lambda,\mu}\left(  x,y\right)  :=\mu
y\left(  x-y\right)  \label{e2}%
\end{align}
Naturally, this map generates a discrete dynamical system - actually a
discrete semidynamical system - in terms of its forward iterates determined by
$n$-fold compositions of the map with itself, denoted as $\Phi_{\lambda,\mu
}^{n}$, or more simply as $\Phi^{n}$, where $n$ is a nonnegative integer. We
shall employ the usual notation and definitions for this discrete system; for
example, the \emph{positive semiorbit }of a point $p\in\mathbb{R}^{2}$, which
we denote as $O_{+}(p)$, is simply defined as
\[
O_{+}(p):=\left\{  \Phi^{n}(p):n\in\mathbb{Z},\,n\geq0\right\}  ,
\]
and all other relevant definitions are standard (\emph{cf}. \cite{guho,
hartman, kathas, palmel, wigbook}). Our model map is clearly real-analytic,
which we denote as usual by $\Phi\in C^{\omega}$ and \emph{a fortiori} smooth,
denoted as $\Phi\in C^{\infty}$.

Assuming the reset and set values, corresponding to the $x$- and
$y$-coordinates, respectively, are normalized so that they may assume the
discrete (logical) values $0$ or $1$, it makes sense to restrict the model map
to the square $I^{2}:=I\times I:=[0,1]\times\lbrack0,1]$ in the plane. In
fact, owing to the obvious symmetry of the circuit with respect to the reset
and set inputs, it is actually natural to further restrict our attention to
the triangular domain%
\begin{equation}
T:=\left\{  (x,y)\in\mathbb{R}^{2}:0\leq x\leq1,\,0\leq y\leq x\right\}  ,
\label{e3}%
\end{equation}
but we shall defer additional discussion of this point until later, except to
note here that%
\begin{equation}
\Phi(0,0)=(1,0),\;\Phi(1,0)=(1,0)\text{ and }\Phi(1,1)=(0,0), \label{e4}%
\end{equation}
which shows that our model is at least logically consistent with the
\textit{R}-\textit{S} flip-flop circuit.

\subsection{Basic properties of the model}

Before embarking on a more thorough dynamical analysis of the iterates of our
simple map model, we shall describe some of its simpler properties. The fixed
points of the map, satisfying $\Phi(x,y)=(x,y)$, are the solutions of the
equations%
\begin{align}
1-x\left[  \lambda\left(  1-x\right)  +y\right]   &  =x\nonumber\\
\mu y\left(  x-y\right)   &  =y, \label{e5}%
\end{align}
from which we readily deduce the property.

\begin{itemize}
\item[(B1)] The points $(1,0)$ and $(1/\lambda,0)$ are fixed points of $\Phi$
for all $\lambda,\mu>0$ (\emph{cf}. (4)), while all other fixed points are in
the complement of the $x$-axis and are determined by the equations%
\[
\left(  1-\lambda\right)  x^{2}+\left(  1+\lambda-\mu^{-1}\right)
x-1=0,\;y=x-\mu^{-1}.
\]
Hence, we have the following: there are no additional fixed points if
$\lambda=1$ and $\mu^{-1}=1+\lambda=2$; there is one more fixed point, $(x,y)$
with
\[
x=\frac{\mu}{\mu(1+\lambda)-1},\;y=\frac{\mu}{\mu(1+\lambda)-1}-\frac{1}{\mu
};
\]
if $\lambda=1$ and $\mu^{-1}\neq1+\lambda=2$; and if $\lambda\neq1$ and
$\left(  1+\lambda-\mu^{-1}\right)  ^{2}+4\left(  1-\lambda\right)  \geq0 $,
there two additional fixed points $(x,y)$ with
\begin{equation}
x=\frac{\left(  \mu^{-1}-\lambda-1\right)  \pm\sqrt{\left(  1+\lambda-\mu
^{-1}\right)  ^{2}+4\left(  1-\lambda\right)  }}{2\left(  1-\lambda\right)
},\;y=x-\frac{1}{\mu}, \label{e6}%
\end{equation}
while if $\left(  1+\lambda-\mu^{-1}\right)  ^{2}+4\left(  1-\lambda\right)
<0$, there are only the two fixed points on the $x$-axis.
\end{itemize}

\noindent The following additional properties of the map follow directly from
its definition.

\begin{itemize}
\item[(B2)] $\Phi$ maps the $x$-axis into itself, and if $0<\lambda\leq4$,
$\Phi$ actually maps the horizontal edge $e_{h}:=\{(x,0):0\leq x\leq1\}$ of
$T$ into itself.

\item[(B3)] $\Phi$ maps the diagonal line $x-y=0$ into the $x$-axis, and maps
the diagonal edge $e_{d}:=\{(x,x):0\leq x\leq1\}$ of $T$ into $e_{h}$ if
$0<\lambda\leq2$.

\item[(B4)] $\Phi$ maps the $y$-axis onto the portion of the line $x=1$ with
$y\leq0$, and maps the line $x=1$, containing the vertical edge $e_{v}%
:=\{(1,y):0\leq y\leq1\}$ of $T$, onto the parabola $y=\mu x(1-x)$ passing
through the origin and the fixed point $(1,0)$.

\item[(B5)] It follows from the derivative (matrix)%
\begin{equation}
\Phi^{\prime}(x,y)=\left(
\begin{array}
[c]{cc}%
\lambda(2x-1)-y & -x\\
\mu y & \mu(x-2y)
\end{array}
\right)  \label{e7}%
\end{equation}
and the inverse function theorem that $\Phi$ is a local $C^{\omega}%
$-diffeomorphism at any point in the complement of the quadratic curve
\[
\lambda\left(  2x-1\right)  \left(  x-2y\right)  +2y^{2}=0,
\]
while in general, the preimage of any point in the plane, denoted as
$\Phi^{-1}\left(  (x,y)\right)  $, is comprised of at most four points.
\end{itemize}

\section{Elementary Dynamics of the Model}

We shall analyze the deeper dynamical aspects of the model map (1)-(2) for
various parameter ranges in the sequel, but first we dispose of some of the
more elementary properties such as a the usual local linear stability analysis
of the fixed points. At this stage, and for the remainder of our
investigation, we shall focus on the restriction of the model map to the
triangle $T$ and assume that

\medskip

\qquad($\mathcal{A}$1)$\qquad\qquad0<\lambda<1<\mu$

\medskip

With the above restriction and assumption, it follows from the preceding
section that our model map has precisely four fixed points: two in $T$; one
near $T$ at $(1/\lambda,0)$ when $\lambda$ is close to unity; and the final
one rather distant from $T$. In the next subsection, we embark on a local
stability analysis of the fixed points of $\Phi$ on or near the triangle $T$.

\subsection{Local analysis of the fixed points}

The local properties of the fixed points of our model map shall be delineated
in a series of lemmas. They all have straightforward proofs that follow
directly from the results in the preceding section and fundamental dynamical
systems theory (as in \cite{guho, kathas, palmel, wigbook}), which are left to
the reader.

\medskip

\noindent\textbf{Lemma 3.1. }\emph{The fixed points of }$\Phi$\emph{\ on the
}$x$\emph{-axis, namely }$(1,0)$\emph{\ and }$(1/\lambda,0)$\emph{, are a
saddle and a source with eigenvalues }(\emph{of }$\Phi^{\prime}(1,0)$
\emph{and} $\Phi^{\prime}(1/\lambda,0)$) $\lambda,\mu$\emph{\ and }%
$2-\lambda,\mu/\lambda$\emph{, respectively. For }$(1,0)$\emph{, the stable
manifold is}%
\[
W^{s}\left(  1,0\right)  =\left\{  (x,0):x<1/\lambda\right\}  ,
\]
\emph{and the linear unstable manifold is}%
\[
W_{\ell}^{u}(1,0)=\left\{  \left(  x,(\lambda-\mu)(x-1)\right)  :x\in
\mathbb{R}\right\}  .
\]

\noindent\textbf{Lemma 3.2. }\emph{The fixed point of }$\Phi$\emph{\ in the
interior of }$T$\emph{, which we denote as }$p_{\ast}=(x_{\ast},y_{\ast}%
)$\emph{, is defined according to }(6) \emph{as }%

\[
x_{\ast}=x_{\ast}(\lambda,\mu)=\frac{\left(  \mu^{-1}-\lambda-1\right)
+\sqrt{\left(  1+\lambda-\mu^{-1}\right)  ^{2}+4\left(  1-\lambda\right)  }%
}{2\left(  1-\lambda\right)  },\;y_{\ast}=y_{\ast}(\lambda,\mu)=x_{\ast}%
-\frac{1}{\mu},
\]
\emph{and has complex conjugate eigenvalues that are roots of the quadratic
equation}%
\[
\sigma^{2}-a\sigma+b=0,
\]
\emph{where}%
\begin{align*}
a  &  =a(\lambda,\mu):=\left(  2\lambda-\mu-1\right)  x_{\ast}+\left(
2-\lambda+\mu^{-1}\right)  ,\;\\
b  &  =b(\lambda,\mu):=\mu\left\{  \left[  \lambda\left(  4\mu^{-1}-1\right)
-2(1+\mu^{-1})\right]  x_{\ast}+2\left[  1-\mu^{-1}\left(  \lambda-\mu
^{-1}\right)  \right]  \right\}  ;
\end{align*}
\emph{namely}%
\begin{align*}
\sigma &  =\sigma(\lambda,\mu)=\frac{1}{2}\left[  a+i\sqrt{4b-a^{2}}\right] \\
\bar{\sigma}  &  =\bar{\sigma}(\lambda,\mu)=\frac{1}{2}\left[  a-i\sqrt
{4b-a^{2}}\right]  .
\end{align*}
\emph{Hence it is a spiral sink or spiral source, respectively, when}%
\[
\left\vert \sigma\right\vert ^{2}=\left\vert \bar{\sigma}\right\vert
^{2}=b<1,
\]
\emph{or}%
\[
\left\vert \sigma\right\vert ^{2}=\left\vert \bar{\sigma}\right\vert
^{2}=b>1.
\]
\emph{Otherwise }(\emph{when }$b=1$) \emph{it has neutral stability.}

\medskip

\noindent\textbf{Lemma 3.3. }\emph{For any fixed }$\lambda$\emph{\ and
variable }$\mu$ \emph{satisfying (}$\mathcal{A}$\emph{1), the coefficient }%
$b$\emph{\ defined above satisfies the following properties}:

\begin{itemize}
\item[(i)] \emph{It is a smooth }$(=C^{\infty})$\emph{, nonnegative function
of }$\mu$\emph{\ for }$\mu>1$\emph{, such that }$db/d\mu>0$ \emph{for every}
$\mu>1$.

\item[(ii)] $b\uparrow\infty$\emph{\ as }$\mu\uparrow\infty.$

\item[(iii)] \emph{There exists a positive }$c(\lambda)$\emph{\ such that
}$b<1$\emph{\ for }$1<\mu<1+c(\lambda).$

\item[(iv)] \emph{In particular, for each }$0<\lambda<1$\emph{, there exists a
unique }$\mu_{h}=\mu_{h}(\lambda)=1+c(\lambda)$\emph{\ such that }$1<\mu_{h}%
$\emph{, }$b(\lambda,\mu_{h})=1$\emph{, }$0<b(\lambda,\mu)<1$\emph{\ for
}$1<\mu<\mu_{h}$\emph{, and }$b(\lambda,\mu)>1$\emph{\ for }$\mu_{h}<\mu
$\emph{.}
\end{itemize}

The case when the interior fixed point is a spiral attractor is shown in Fig.
2. By fixing $\lambda$ at a value near one, say $\lambda=0.99$, and then
varying $\mu$ over a range from $4$ to $5$, we obtain a very rich array of
dynamics as described in what follows.

\section{Oscillation and Hopf Bifurcation}

In order to achieve a reasonable amount of focus - given the wide range of
possible model map parameters - we shall narrow our range of investigation by
adhering to the following additional assumption in the sequel:

\medskip

\qquad($\mathcal{A}$2)$\qquad\qquad\lambda=0.99=\frac{99}{100}.$

\medskip

\noindent Then assuming ($\mathcal{A}$1) and ($\mathcal{A}$2), our map
$\Phi=\Phi_{\mu}$ satisfies all the properties delineated above, and depends
only on the single parameter $\mu\in(1,\infty)$. In particular, it follows
directly from Lemma 3.2 that%
\begin{equation}
x_{\ast}(\mu):=x_{\ast}(.99,\mu)=50\left\{  \left(  \mu^{-1}-1.99\right)
+\sqrt{\left(  1.99-\mu^{-1}\right)  ^{2}+0.04}\right\}  ,\;y_{\ast}%
(\mu):=y_{\ast}(.99,\mu)=x_{\ast}-\frac{1}{\mu}, \label{e8}%
\end{equation}
and%
\begin{align}
a(\mu)  &  :=a(.99,\mu)=\left(  .98-\mu\right)  x_{\ast}+\left(  1.01+\mu
^{-1}\right)  ,\nonumber\\
b(\mu)  &  :=b(.99,\mu)=\mu\left\{  \left[  (.99)\left(  4\mu^{-1}-1\right)
-2(1+\mu^{-1})\right]  x_{\ast}+2\left[  1-\mu^{-1}\left(  .99-\mu
^{-1}\right)  \right]  \right\}  .\; \label{e9}%
\end{align}
It is then straightforward to compute in the notation of Lemma 3.3 that%
\begin{equation}
c:=c(.99)\cong3.5438,\;\mu_{h}:\cong\mu_{h}(.99)\cong4.5438,\;x_{\ast}(\mu
_{h})\cong0.5632,\;y_{\ast}(\mu_{h})\cong0.3431. \label{e10}%
\end{equation}

In order to study the behavior of the map in a neighborhood of the fixed point
$p_{\ast}=(x_{\ast},y_{\ast})$, it is convenient to translate the coordinates
and map to the origin by defining%
\begin{align}
\hat{\Phi}  &  :=\hat{\Phi}_{\mu}(\xi,\eta):=\Phi_{\mu}(\xi+x_{\ast}%
,\eta+y_{\ast})-\Phi_{\mu}(x_{\ast},y_{\ast})\nonumber\\
&  =\Phi_{\mu}(\xi+x_{\ast},\eta+y_{\ast})-(x_{\ast},y_{\ast}). \label{e11}%
\end{align}
It is easy to compute that%
\begin{equation}
\hat{\Phi}:=\hat{\Phi}_{\mu}(\xi,\eta)=\left(  \hat{\varphi}_{\mu}(\xi
,\eta),\hat{\psi}_{\mu}(\xi,\eta)\right)  , \label{e12}%
\end{equation}
where%
\begin{align}
\hat{\varphi}_{\mu}(\xi,\eta)  &  :=\left[  (0.98)x_{\ast}(\mu)+\mu
^{-1}-0.99\right]  \xi-x_{\ast}(\mu)\eta+\xi\left[  (0.99)\xi-\eta\right]
,\nonumber\\
\hat{\psi}_{\mu}(\xi,\eta)  &  :=\left[  \mu x_{\ast}(\mu)-1\right]
\xi-\left[  2-\mu x_{\ast}(\mu)\right]  \eta+\mu\eta\left(  \xi-\eta\right)  .
\label{e13}%
\end{align}

\subsection{Invariant curve and Hopf bifurcation}

Now with this simple quadratic representation of the model map with respect to
the fixed point $p_{\ast}$ interior to the triangle $T$ is a straightforward
matter to describe the bifurcation behavior and oscillatory properties. In
particular, we have the following result.

\medskip

\noindent\textbf{Theorem 4.1. }\emph{The discrete semidynamical system
associated to the map }$\Phi_{\mu}$\emph{\ }$($\emph{or }$\hat{\Phi}_{\mu}%
)$\emph{\ has a Hopf bifurcation at the fixed point }$p_{\ast}$\emph{\ when
}$\mu=\mu_{h}$\emph{. More specifically, }$p_{\ast}$\emph{\ is a spiral sink
}$($\emph{source}$)$\emph{\ for }$1<\mu<\mu_{h}$\emph{\ }$(\mu_{h}<\mu
)$\emph{, and }$p_{\ast}$\emph{\ is neutrally stable for }$\mu=\mu_{h}%
$\emph{\ with }$\Phi_{\mu_{h}}^{\prime}(p_{\ast})$\emph{\ having complex
conjugate eigenvalues on the unit circle }$S^{1}$\emph{\ in the complex plane
}$\mathbb{C}$\emph{\ }$(``=$\textquotedblright\ $\mathbb{R}^{2})$\emph{.
Furthermore, for sufficiently small }$\nu:=\mu-\mu_{h}>0$\emph{, say }%
$0<\nu<\epsilon$\emph{, there exists a unique }$\Phi_{\mu}$\emph{-invariant,
smooth Jordan curve }$\Gamma_{\nu}$\emph{\ }$($\emph{i.e. with }$\Phi_{\mu
}(\Gamma_{\nu})=\Gamma_{\nu})$\emph{\ enclosing }$p_{\ast}$\emph{\ in its
interior, }$\mathcal{I}(\Gamma_{\nu}),$\emph{satisfying the following
properties:}

\begin{itemize}
\item[(i)] \emph{There is a }$0<\epsilon_{\ast}\leq\epsilon$\emph{\ for which
}$0<\nu<\epsilon_{\ast}$\emph{\ implies that for every point }$p\in
\mathcal{I}(\Gamma_{\nu})\smallsetminus\{p_{\ast}\}$\emph{\ the iterates
}$\Phi_{\mu}^{n}(p)$\emph{\ spiral around }$p_{\ast}$\emph{\ in a
counterclockwise manner and approach }$\Gamma_{\nu}$\emph{; in particular, the
distance between these iterates and the curve, denoted }$\Delta(\Phi_{\mu
_{h}+\nu}^{n}(p),\Gamma_{\nu})$\emph{, converges monotonically to zero as
}$n\rightarrow\infty$\emph{.}

\item[(ii)] \emph{With }$\epsilon_{\ast}$\emph{\ and }$\epsilon$\emph{\ as in
}$(i)$\emph{, }$\Gamma_{\nu}$\emph{\ is a local attractor, in that all
positive semiorbits originating in some open neighborhood of this curve
converge to }$\Gamma_{\nu}$\emph{.}

\item[(iii)] \emph{The dynamical system on }$\Gamma_{\nu}$\emph{\ induced by
the restriction of the map }$\Phi_{\mu}$\emph{\ is either ergodic }%
$($\emph{with dense orbits}$)$\emph{\ or has periodic orbits }$($\emph{or
cycles}$)$\emph{\ according as the rotation number is irrational or rational,
respectively, including an }$11$\emph{-cycle as }$\mu$\emph{\ just exceeds the
bifurcation value }$\mu_{h}$\emph{, the collection }$\{\Gamma_{\nu}%
:0<\nu<\epsilon_{\ast}\}$\emph{\ includes }$m$\emph{-cycles for infinitely
many }$m\in\mathbb{N}$\emph{, where }$\mathbb{N}$\emph{\ comprises the natural
numbers.}

\item[(iv)] \emph{Under the same conditions as in }$(i)$\emph{, }%
$\Delta(p_{\ast},\Gamma_{\nu})=O(\nu).$
\end{itemize}

\smallskip

\noindent\emph{Proof. }Properties (i) and (ii) follow from a direct
application to the map (11) of the Hopf bifurcation theorem for discrete
dynamical systems (see \emph{e.g.} \cite{DMP, MM, wan, WXH} and also
\cite{guho, IJ, kathas, wigbook}), or one can obtain the same results via a
straightforward modification of the main theorem of Champanerkar \& Blackmore
\cite{CB}. In fact, the latter approach actually shows that the invariant
curve is analytic in its variables and parameter $\nu$. To prove (iii) and
(iv) requires a deeper analysis of the invariant curves, which we shall merely
outline in the interest of brevity (\emph{cf}. \cite{guho} and Lanford's
version of Ruelle's proof in \cite{MM}).

The curve $\Gamma_{\nu}$ can be parametrized in polar form as%
\begin{equation}
\Gamma_{\nu}:\xi=\xi(\theta;\mu):=\rho(\theta;\mu)\cos\theta,\quad\eta
=\eta(\theta;\mu):=\rho(\theta;\mu)\sin\theta, \label{e14}%
\end{equation}
where $\theta$ is the usual polar angle about the point $p_{\ast}$. Then
$\Phi_{\mu}$-invariance requires that
\begin{equation}
\hat{\Phi}_{\mu}\left(  \xi(\theta;\mu),\eta(\theta;\mu)\right)  =\left(
\xi(\Theta;\mu),\eta(\Theta;\mu)\right)  , \label{e15}%
\end{equation}
where $\Theta$ represents the angular rotational action of the map defined as%
\begin{equation}
\Theta=\Theta\left(  \xi,\eta;\nu\right)  :=\tan^{-1}\left(  \frac{\hat{\psi
}_{\mu}(\xi,\eta)}{\hat{\varphi}_{\mu}(\xi,\eta)}\right)  . \label{e16}%
\end{equation}
We note that is easy to verify that we must take the branch of the
arctangent that takes on values between $\theta+\pi/2$ and $\theta+3\pi/2.$

It is convenient to introduce the following more compact notation for the map
as expressed in terms of its coordinate functions in (13):%
\begin{align}
\hat{\varphi}_{\mu}(\xi,\eta)  &  :=-\alpha(\nu)\xi-\beta(\nu)\eta+\xi\left[
(0.99)\xi-\eta\right]  ,\nonumber\\
\hat{\psi}_{\mu}(\xi,\eta)  &  :=\gamma(\nu)\xi+\delta(\nu)\eta+\mu\eta\left(
\xi-\eta\right)  , \label{e17}%
\end{align}
where the parameter dependent, positive coefficients $\alpha(\nu)$, $\beta
(\nu)$, $\gamma(\nu)$ and $\delta(\nu)$ are defined in the obvious way
according to (13). In polar coordinates with respect to $p_{\ast}$, the
functions defined in (17) have the form%
\begin{align}
\hat{\varphi}_{\mu}(r,\theta)  &  :=r\left\{  -\alpha(\nu)\cos\theta-\beta
(\nu)\sin\theta+r\cos\theta\left[  (0.99)\cos\theta-\sin\theta\right]
\right\}  ,\nonumber\\
\hat{\psi}_{\mu}(r,\theta)  &  :=r\left[  \gamma(\nu)\cos\theta+\delta
(\nu)\sin\theta+\mu r\sin\theta\left(  \cos\theta-\sin\theta\right)  \right]
. \label{e18}%
\end{align}
Moreover, in the context of these polar coordinates, our map can be rewritten
in the form%
\begin{equation}
\hat{\Phi}_{\mu}\left(  r,\theta\right)  :=\left(  R(r,\theta;\mu
),\Theta(r,\theta;\mu)\right)  , \label{e19}%
\end{equation}
and we compute that%
\begin{align*}
R^{2}  &  =\hat{\varphi}_{\mu}^{2}+\hat{\psi}_{\mu}^{2}=r^{2}\left\{  \left(
\alpha^{2}+\gamma^{2}\right)  \cos^{2}\theta+\left(  \alpha\beta+\gamma
\delta\right)  \cos2\theta+\left(  \beta^{2}+\delta^{2}\right)  \sin^{2}%
\theta-\right. \\
&  r\left(  0.99\cos\theta-\sin\theta\right)  \left[  2\alpha\cos^{2}%
\theta+\beta\cos2\theta-r\cos^{2}\theta\left(  0.99\cos\theta-\sin
\theta\right)  \right]  +\\
&  \left.  \mu r\left(  \cos\theta-\sin\theta\right)  \left[  \gamma
\cos2\theta+2\delta\sin^{2}\theta+r\sin^{2}\theta\left(  \cos\theta-\sin
\theta\right)  \right]  \right\}  ,
\end{align*}
which implies that%
\begin{align}
R\left(  r,\theta;\mu\right)   &  :=rU\left(  r,\theta;\mu\right)  =r\left\{
\left(  \alpha^{2}+\gamma^{2}\right)  \cos^{2}\theta+\left(  \alpha
\beta+\gamma\delta\right)  \cos2\theta+\left(  \beta^{2}+\delta^{2}\right)
\sin^{2}\theta-\right. \nonumber\\
&  r\left(  0.99\cos\theta-\sin\theta\right)  \left[  2\alpha\cos^{2}%
\theta+\beta\cos2\theta-r\cos^{2}\theta\left(  0.99\cos\theta-\sin
\theta\right)  \right]  +\label{e20}\\
&  \left.  \mu r\left(  \cos\theta-\sin\theta\right)  \left[  \gamma
\cos2\theta+2\delta\sin^{2}\theta+r\sin^{2}\theta\left(  \cos\theta-\sin
\theta\right)  \right]  \right\}  ^{1/2}\nonumber
\end{align}

Now it follows from the $\hat{\Phi}_{\mu}$-invariance of $\Gamma_{\nu}$,
manifested by (19), and (14)-(20) that the radius function $\rho$ must satisfy%
\begin{equation}
\rho=\rho\left(  \theta;\mu\right)  =U\left(  \rho,\theta;\mu\right)
^{-1}\rho\left(  \tan^{-1}\left(  \frac{\left[  \gamma(\nu)\cos\theta
+\delta(\nu)\sin\theta+\mu\rho\sin\theta\left(  \cos\theta-\sin\theta\right)
\right]  }{\left\{  -\alpha(\nu)\cos\theta-\beta(\nu)\sin\theta+\rho\cos
\theta\left[  (0.99)\cos\theta-\sin\theta\right]  \right\}  }\right)  \right)
. \label{e21}%
\end{equation}
This equation expresses the fact that $\rho$ is a fixed point of the operator
on the right-hand side, and can be used to approximate this radius function to
any desired degree of accuracy. As we noted above, $\rho$ is analytic in
$(\theta,\mu)$, so the approximation can be effected by assuming a power
series representation in $\theta$ or $\mu$, whereupon substitution in (21)
would provide a means for recursive determination of the series coefficients.
But this turns out to be a rather laborious, albeit straightforward, process.
It is actually more efficient in this case to find (global in $\theta$)
approximations via Picard iteration. For example, if we take $\rho_{0}=\nu$,
the next approximation yields%
\[
\rho_{1}=\nu U\left(  \nu,\theta;\mu\right)  ^{-1},
\]
which owing to the easily verified rather rapid geometric convergence of the
iterates, gives quite a good approximation of the (fixed point) solution - one
that is, for example, sufficient to verify property (iv). Then a closer
examination of the accuracy of the remaining successive approximations, with
special attention to the angular aspect of the restriction of $\hat{\Phi}%
_{\mu}$ to $\Gamma_{\nu}$ embodied in%
\[
\Theta=\tan^{-1}\left(  \frac{\left[  \gamma(\nu)\cos\theta+\delta(\nu
)\sin\theta+\mu\rho\sin\theta\left(  \cos\theta-\sin\theta\right)  \right]
}{\left\{  -\alpha(\nu)\cos\theta-\beta(\nu)\sin\theta+\rho\cos\theta\left[
(0.99)\cos\theta-\sin\theta\right]  \right\}  }\right)  ,
\]
together with some fundamental results on rotation numbers such as given in
Hartman \cite{hartman}, makes it possible to verify (iii); thereby completing
the proof. $\blacksquare$

\subsection{Cascading Hopf doubling bifurcations}

If one continues further along the lines of analysis of the behavior of the
map $\hat{\Phi}_{\mu}$ in the proof of Theorem 4.1, a much more intricate
sequence of bifurcations and dynamical properties emerges, which we shall just
sketch here. We begin by keeping close tabs on the stability of the invariant
curve $\Gamma_{\nu}$ as $\nu$ increases. A careful analysis of this locally
attracting curve and the map, which we leave to the reader, reveals that there
is a small $\nu_{1}>0$ beyond which the curve becomes locally repelling, and
for which the main theorem of \cite{CB} applies. Accordingly a pair of new
locally attracting, smooth Jordan curves emerge from a pitchfork bifurcation
of $\Gamma_{\nu}$ - one interior to $\Gamma_{\nu}$, which we denote as
$\Gamma_{\nu}^{(0)}$, and the other exterior to $\Gamma_{\nu}$, denoted as
$\Gamma_{\nu}^{(1)}$, such that $\Gamma_{\nu}^{(0)}\cup\Gamma_{\nu}^{(1)}$ is
$\hat{\Phi}_{\mu}$-invariant, with $\hat{\Phi}_{\mu}(\Gamma_{\nu}%
^{(0)})=\Gamma_{\nu}^{(1)}$ and $\hat{\Phi}_{\mu}(\Gamma_{\nu}^{(1)}%
)=\Gamma_{\nu}^{(0)}$. In effect then, $\{\Gamma_{\nu}^{(0)},\Gamma_{\nu
}^{(1)}\}$ is a 2-cycle of sets. Thus we have a doubling bifurcation for
one-dimensional closed smooth manifolds analogous to the beginning of a period
doubling cascade for points (zero-dimensional manifolds) observed in such
one-dimensional discrete dynamical systems as that of the logistic map.

One can actually show that this analog is complete, in that there is an
infinite sequence of such stability shifting, smooth, invariant, Jordan curve
doubling bifurcations that converge to an extremely complicated chaotic state.
More specifically, there is a $\nu_{2}>\nu_{1}$ such that across this
parameter value, each member of the pair $\Gamma_{\nu}^{(0)}$, $\Gamma_{\nu
}^{(1)}$ becomes locally repelling, and gives birth - via pitchfork
bifurcation (\emph{cf}. \cite{CB}) - to a pair of locally attracting, smooth
Jordan curves, $\Gamma_{\nu}^{(0,0)},\Gamma_{\nu}^{(0,1)}$ and $\Gamma_{\nu
}^{(1,0)},\Gamma_{\nu}^{(1,1)}$, respectively. Furthermore, $\Gamma_{\nu
}^{(0,0)}$ ($\Gamma_{\nu}^{(0,1)}$) is in the interior (exterior) of
$\Gamma_{\nu}^{(0)}$ and $\Gamma_{\nu}^{(1,0)}$ ($\Gamma_{\nu}^{(1,1)}$) is in
the interior (exterior) of $\Gamma_{\nu}^{(1)}$, $\Gamma_{\nu}^{(0,0)}%
\cup\Gamma_{\nu}^{(0,1)}\cup\Gamma_{\nu}^{(1,0)}\cup\Gamma_{\nu}^{(1,1)}$ is
$\hat{\Phi}_{\mu}$-invariant, with these four curves forming the 4-cycle
$\Gamma_{\nu}^{(0,0)}\rightarrow\Gamma_{\nu}^{(1,0)}\rightarrow\Gamma_{\nu
}^{(0,1)}\rightarrow$ $\Gamma_{\nu}^{(1,1)}\rightarrow\Gamma_{\nu}^{(0,0)}$.
This process continues \emph{ad infinitum }to generate a bounded monotone
increasing sequence of parameter values $\nu_{1}<\nu_{2}<\nu_{3}<\cdots$, with
$\nu_{n}\rightarrow\nu_{\infty}<3$, with $\Phi_{\mu_{h}+\nu_{\infty}}$
exhibiting chaotic dynamics. Among the consequences of this cascade of
bifurcations is that for $\mu\geq\mu_{h}+\nu_{\infty}$ there exists a closed,
smooth, $\hat{\Phi}_{\mu}$-invariant curvilinear annulus $A$ that encloses the
fixed point $p_{\ast}$ and contains an infinite number (with cardinality of
the continuum) of smooth Jordan curves that can be partitioned into
$n$-cycles, with $n$ ranging over the nonnegative integers, and naturally this
annulus contains very intricate dynamics. We summarize this in the next result
- illustrated in Fig. 3 for the Hopf bifurcation and Fig. 4 for the multiring
configuration - whose proof we shall leave to the reader for now, although we
plan to prove it in a more general form in a forthcoming paper. It is helpful
to observe that it describes the discrete dynamical behavior embodied in the
following paradigm represented in polar coordinates, with angular coordinate
function similar to the circular map of Arnold \cite{VIA}:%
\[
(r,\theta)\rightarrow\left(  R,\Theta\right)  ,
\]
where%
\begin{align*}
R &  :=\left(  \nu+1\right)  r\left(  1-r\right)  ,\\
\Theta &  :=\theta+\frac{2\pi}{a+\nu}\left(  1+kr\sin\theta\right)
\;(\mathrm{mod\,}2\pi),
\end{align*}
and $a$ and $k$ are positive numbers.

\medskip

\noindent\textbf{Theorem 4.2. }\emph{Let }$\hat{\Phi}_{\nu}$ \emph{be defined
as }$\hat{\Phi}_{\mu_{h}+\nu}$. \emph{There exists an increasing sequence
}$0=\nu_{0}<\nu_{1}<\nu_{2}<\cdots\rightarrow\nu_{\infty}<3$\emph{\ such that,
in addition to the }$\hat{\Phi}_{\nu}$-\emph{invariant smooth Jordan curve
}$\Gamma_{\nu}$,\emph{\ which exists for all }$\nu>\nu_{0}$,\emph{\ and is
locally attracting }$(\emph{repelling})$\emph{\ for }$\nu_{0}<\nu<\nu_{1}%
$\emph{\ }$(\nu>\nu_{1})$,\emph{\ for each }$m\in\mathbb{N}$\emph{\ and }%
$\nu>\nu_{m}$ \emph{there is a }$2^{m}$-\emph{\ cycle }$($\emph{of sets that
are smooth Jordan curves}$)$\emph{\ of }$\hat{\Phi}_{\nu},$%
\[
\mathcal{Z}_{m}:=\left\{  \Gamma_{\nu}^{i_{m}}:i_{m}\in\{0,1\}^{m}\right\}
\]
\emph{\ which is created from }$\mathcal{Z}_{m-1}$ \emph{via pitchfork
bifurcation, and is locally attracting }$($\emph{repelling}$)$\emph{\ for
}$\nu_{m}<\nu<\nu_{m+1}$ $\emph{(}\nu>\nu_{m+1})$.\emph{\ Furthermore, }%
$\hat{\Phi}_{\nu}^{2^{m}}$\emph{\ restricted to any of the curves }%
$\Gamma_{\nu}^{i_{m}}$\emph{\ is either ergodic or has periodic orbits
according as the rotation number }$($\emph{which varies continuously with
}$\nu)$\emph{\ is irrational or rational, respectively. Consequently, }%
$\hat{\Phi}_{\nu}$ \emph{has periodic orbits of arbitrarily large period for
infinitely many }$\nu$\emph{\ in a small neighborhood of \ }$\nu_{\infty}$,
\emph{where it can also be shown to have chaotic orbits if }$\nu>\nu_{\infty}%
$\emph{. Furthermore, for }$\nu\geq\nu_{\infty}$ \emph{there is a closed,
smooth, }$\hat{\Phi}_{\nu}$-\emph{invariant annulus }$A$ \emph{\ enclosing the
fixed point }$p_{\ast}$, \emph{which contains }$\Gamma_{\nu}$\emph{, is
locally attracting and is the minimal invariant set containing all the cycles
}$\mathcal{Z}_{m}$.

\medskip

\noindent The invariant, locally attracting annulus $A$ may not precisely
qualify as a strange attractor, yet the intricacies of the dynamics it
contains - including cycles of arbitrarily large period - deserves a special
name such as a \emph{pseudo-strange attractor}.

\begin{figure}[thb]
\centering
\includegraphics[width=4.5in]{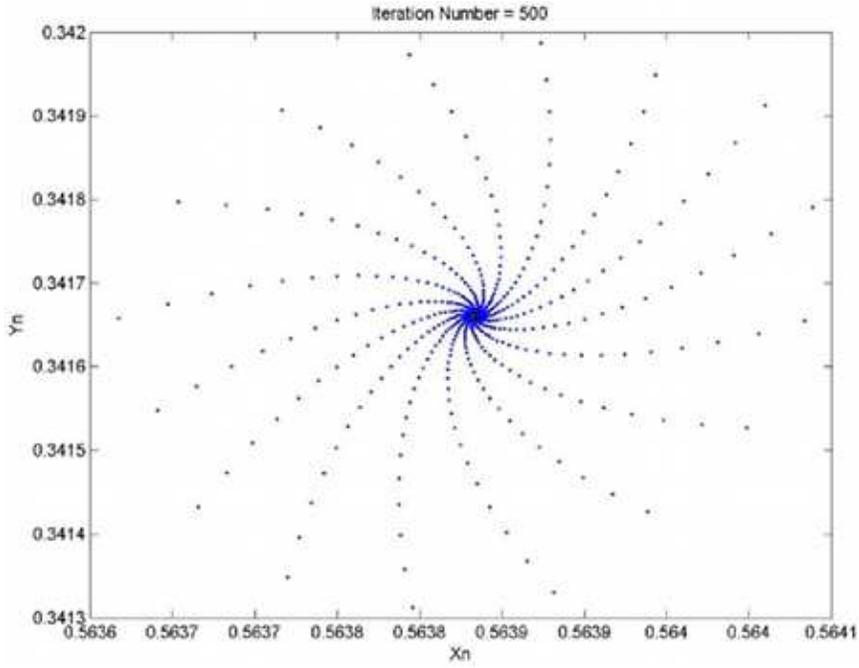}\caption{Spiral attractor at
interior fixed point ($\cong(0.5639, 0.3417)$) for $\lambda=0.99,\, \mu=4.5$.
The initial point is $(x,y)=(0.564,0.342)$ and the eigenvalues are $\xi
\cong-0.3763\pm0.9171i$.}%
\label{attractor}%
\end{figure}

\begin{figure}[hb]
\centering
\includegraphics[width=4.5in]{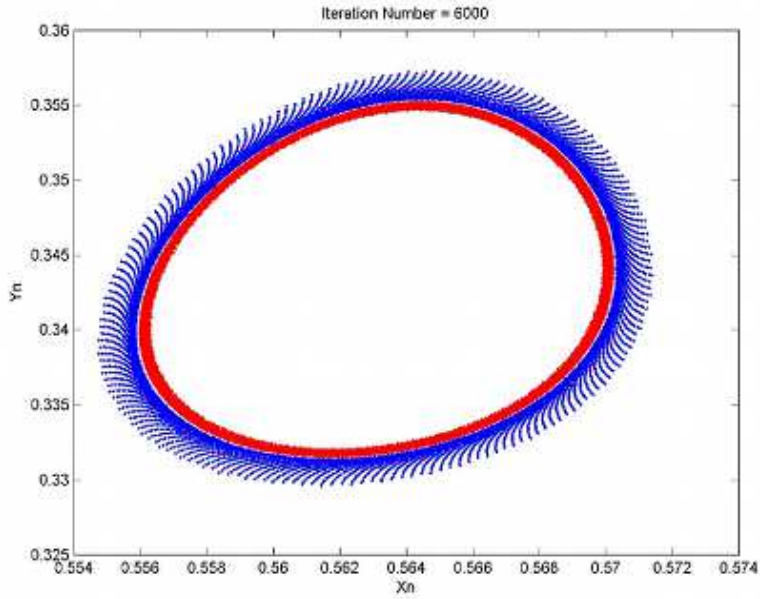}\caption{Hopf bifurcation at
interior fixed point ($\cong(0.5639, 0.3417)$) for $\lambda=0.99,\,
\mu=4.5449$, with eigenvalues $\xi\cong-0.3889\pm0.9215i$. The initial points
are $(x,y)=(0.555,0.340)$ and $(x,y)=(0.558,0.34)$.}%
\label{Hopf}%
\end{figure}

\begin{figure}[bht]
\centering
\includegraphics[width=4.5in]{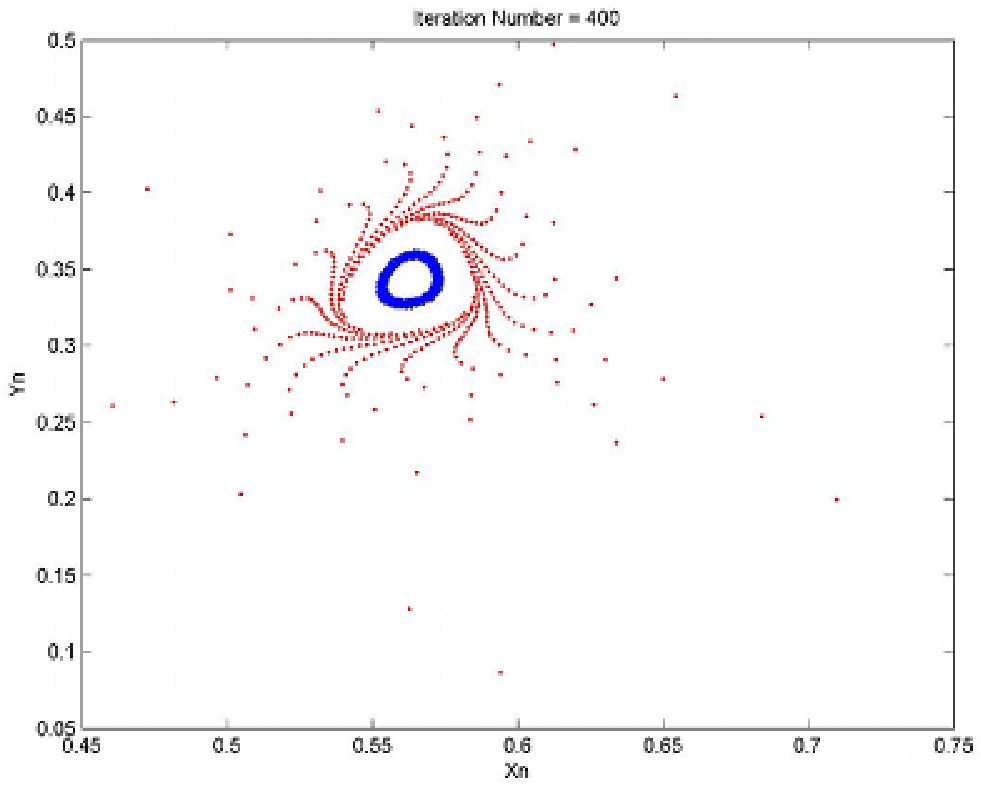}\caption{Multiring structure
around interior fixed point ($\cong(0.5639, 0.3417)$) for $\lambda=0.99,\,
\mu=4.55$, with eigenvalues $\xi\cong-0.3903\pm0.922i$. The initial points are
$(x,y)=(0.555,0.338)$ and $(x,y)=(0.43,0.38)$.}%
\label{multiring}%
\end{figure}

\section{Chaotic Dynamics and Instability}

As pointed out in Theorem 4.2, our model exhibits chaotic dynamics at the
limit of the cascade of doubling bifurcations described therein. The proof of
this \textquotedblleft limiting form\textquotedblright\ of chaos turns out to
be rather subtle and difficult, so we shall not go into it here. Instead, we
shall prove the existence of chaotic regimes for higher parameter values using
a fairly simple geometric argument based upon demonstrating that a
sufficiently high power (iterate) of our model map exhibits Smale
horseshoe-like behavior (see \emph{e.g.} \cite{guho, kathas, palmel,
wigbook}), as illustrated in Fig. 5, and also heteroclinic cycles with
transverse intersections (\emph{cf}. \cite{db, palmel, wigbook}). We keep the
value of $\lambda=0.99$, and set $\mu=5$, and summarize our main findings on
chaos in the following result.

\medskip

\noindent\textbf{Theorem 5.1. }\emph{Let }$\Phi:=\Phi_{0.99,5}$. \emph{Then
there exist an }$N\in\mathbb{N}$ \emph{and a region }$Q$\emph{\ homeomorphic
with a square that intersects the triangular domain }$T$, \emph{and has the
following properties depicted in Fig. }$5$:

\begin{itemize}
\item[(a)] $F:=\Phi^{N}$\emph{ maps }$Q$\emph{ onto a pill shaped region,
which contains the fixed point }$(1,0)$\emph{ in its interior, lies along the
vertical edge }$e_{v}$\emph{ of }$T$\emph{, has maximum }$y$\emph{ value of
}$0.22$\emph{, and is of sufficient height to include the line segment from
}$(1,0)$\emph{ to }$(1,0.2)$\emph{ in its interior.}

\item[(b)] $G:=\Phi^{N+1}$\emph{ maps }$Q$\emph{ along the curve described
above in (B4) into a thin curved homeomorph }$\Phi(F(Q))$\emph{ of a square
that just crosses the diagonal edge }$e_{d}$\emph{ of }$T$\emph{ in the
interior of the unit square.}

\item[(c)] $H:=\Phi^{N+2}=\Phi\circ G$\emph{ maps }$Q$\emph{ into a stretched
and folded homeomorph of }$Q$\emph{ that (transversely) intersects }$Q$\emph{
in a disjoint pair of curvilinear rectangles in a horseshoe-like fashion.}
\end{itemize}

\noindent\emph{Consequently, }$Q$ \emph{contains a compact invariant set
}$\Lambda$ \emph{on which }$H$ \emph{is topologically conjugate to a shift map
}$($\emph{or equivalently, }$\Phi$ \emph{\ is topologically conjugate to a
subshift}$),$\emph{\ which implies that }$\Phi$ \emph{generates chaotic
dynamics in }$Q\cap\Lambda$\emph{, including a dense orbit and cycles of
arbitrarily large period.}

\medskip

\noindent\emph{Proof}. We begin our proof with a disk $D$ with an elliptical
boundary having its center at $(0.95,0.1)$, semimajor axis length equal to
$0.13$ and semiminor axis length equal to $0.075$. Then it follows from the
properties of the model map delineated in (B1)-(B5) above that there exist a
sufficiently large positive integer $N$ and a diffeomorph of $D$, which we
denote as $Q$, such that $\Phi^{-N}(D)=Q$, where $Q$ contains the horizontal
edge $e_{h}$ of \ $T$ as shown in Fig. 5(a). Observe here that $\Phi$
restricted to $e_{h}$ is a doubling map taking the left vertex into the right
vertex, which is symmetric with respect to $x=1/2$, so the inverse notation in
this definition must be viewed in the set theoretical preimage context.
Nevertheless, $Q$ satisfies $F(Q):=\Phi^{N}(Q)=D$.

Next, we consider the image of $D$ under $\Phi$; namely, $\Phi(D)=\Phi\left(
F(Q)\right)  =G(Q)$, which is illustrated in Fig. 5(b). To see that this is an
accurate depiction of the image, first observe that for the point $(1,0.2)$
lying in the interior of $D$ near its highest point, we have $\Phi
(1,0.2)=(0.8,0.8)$, which lies on the diagonal edge $e_{d}$ of $T$. Moreover,
the fixed point $(1,0)$ is also an interior point of $D$, but one that lies
near its lowest point. Accordingly in virtue of the properties of $\Phi$ - in
particular (B4) - the shape of $\Phi(D)$ must be that of the thickened version
of part of the parabolic curve described in (B4) and must slightly overlap
$e_{d}$ around the point $(0.8,0.8)$, just as depicted in Fig. 5(a).

In order to obtain the desired horseshoe-like behavior, another application of
$\Phi$ is required; that is, we need to describe $\Phi^{2}(D)=\Phi\left(
G(Q)\right)  =H(Q)$. Taking into account the overall definition of the map,
property (B3), and the fact that $\Phi(0.8,0.8)=(0.2016,0)$, it is clear that
Fig. 4(c) is a rather accurate rendering of the region $H(Q)$, which
intersects $Q$ in a fairly typical horseshoe type set comprised of two
approximately rectangular components. With this geometric representation of
$Q$ and its image $H(Q)$ in hand, we can describe the key fractal component
$\Lambda$ of nonwandering set $\Omega$ of $H$ (and \emph{a fortiori}, $\Phi$)
and the topological conjugacy of the restriction of $H$ on $\Lambda$, which we
denote as $h$, to a shift map on doubly-infinite binary sequences in
essentially the usual way (\emph{cf}. \cite{guho, kathas, palmel, wigbook}),
modulo a minor alteration necessitated by the fact that $\Phi$ fails to
injective on some subsets of $T$.

It remains to describe the alteration and the final steps in defining
$\Lambda$ and establishing the topological conjugacy between $h:\Lambda
\rightarrow\Lambda$ and the shift map $\sigma:2^{\mathbb{Z}}\rightarrow
2^{\mathbb{Z}}$, where $2^{\mathbb{Z}}:=\{0,1\}^{\mathbb{Z}}$ - the space of
all doubly-infinite binary sequences $\ldots a_{-2}a_{-1}a_{0}a_{1}a_{2}%
\ldots$with $a_{i}=0$ or $1$. To this end, we define $C_{\ell}$ and $C_{r}$ to
be the left and right components of $H(Q)\cap Q$, respectively, as shown in
Fig. 5(c), and observe that $C_{\ell}\subset Q_{l}:=Q\cap\{(x,y):x<2/5\}$ and
$C_{r}\subset Q_{r}^{+}:=\left(  Q\cap\{(x,y):x>3/5\}\right)  \cup D$. It
follows readily from the definition and properties of $\Phi$ delineated in
Section 2 that by making $Q$ more slender and $N$, larger, if necessary,
$\Phi$ maps $Q_{\ell}$ and $Q_{r}^{+}$ diffeomorphically onto their images,
with $\Phi\left(  Q_{\ell}\right)  ,\Phi\left(  Q\cap\{(x,y):x>3/5\}\right)
\subset$ $\{(x,y):x>3/5\}$. In addition, possibly after another shrinking of
$Q$ and increase of $N$, we may assume that the iterated sets $\{C_{\ell}%
,\Phi^{m}\left(  C_{\ell}\right)  \cap Q_{r}^{+}:m\in\mathbb{N}\}$ are
pairwise disjoint. We denote the restriction of $\Phi$ to $Q_{\ell}$ and
$Q_{r}$ by $\Phi_{\ell}$ and $\Phi_{r}$, respectively. Furthermore, $C_{\ell}$
is the diffeomorphic image of a disjoint set $C_{\ell}^{-1}$ ($=\Phi_{\ell
}^{-1}(C_{\ell})$) under $\Phi$, which intersects the edge $e_{d}$, and
$C_{\ell}^{-1}$ is, in turn, the diffeomorphic image under $\Phi_{r}$ of a set
$C_{\ell}^{-2}$ ($=\Phi_{r}^{-1}(C_{\ell}^{-1})$) contained in $D=\Phi^{N}(Q)$.

We can now define a unique inverse on $\Lambda\subset C_{\ell}\cup C_{r}$ - a
set to be defined in this last phase of our proof. As $H^{-1}=\Phi^{-1}%
\circ\Phi^{-1}\circ\cdots\circ\Phi^{-1}$ ($N$ factors), it remains to select
the proper branch of each of the factors in this composition so as to
prescribe $H^{-1}$ unambiguously. It is easy to see that this is accomplished
as follows: set%
\[
\Phi_{\ast}^{-1}\left(  x,y\right)  :=\left\{
\begin{array}
[c]{cc}%
\Phi_{\ell}^{-1}(x,y), & \left(  x,y\right)  \in C_{\ell}\cup\Phi_{\ell
}\left(  C_{\ell}\right) \\
\Phi_{r}^{-1}(x,y) & (x,y)\in C_{\ell}^{-1}\cup\left(  Q_{r}^{+}%
\smallsetminus\Phi_{\ell}\left(  C_{\ell}\right)  \right)
\end{array}
\right.  ,
\]
and%
\[
H^{-1}:=\Phi_{\ast}^{-1}\circ\Phi_{\ast}^{-1}\circ\cdots\circ\Phi_{\ast}%
^{-1}\quad(N\text{ factors}),
\]
and define
\[
\Lambda:=Q\cap\left[
%TCIMACRO{\dbigcap \nolimits_{m\in\mathbb{Z}}}%
%BeginExpansion
{\displaystyle\bigcap\nolimits_{m\in\mathbb{Z}}}
%EndExpansion
H^{m}\left(  H(Q)\cap Q\right)  \right]  .
\]
It is easy to verify that $\Lambda$ is a compact, $H$-invariant set, which is
homeomorphic with the cartesian product of a pair of 2-component Cantor sets
that , in turn, is homeomorphic with $2^{\mathbb{Z}}$ employing the standard
topologies. Then the restriction $h:=H_{\mid\Lambda}$ can be shown to be
topologically conjugate to the shift map in the usual way (such as in
\cite{guho, palmel, wigbook}), thereby completing the proof. $\blacksquare$

\medskip

We note here that the existence of chaotic regimes described in Theorem 5.1
could also have been demonstrated following a more detailed analysis of the
iterates along the lines of the above proof - revealing both transverse
homoclinic points of periodic points and transverse heteroclinic points of
branches of heteroclinic cycles of periodic points, both of which imply the
existence of chaos (as shown or indicated in \cite{db, guho, kathas, palmel,
wigbook}). The chaotic case (with its characteristic splattering effect) is
depicted in Fig. 6, which shows the iterates corresponding to three initial
points selected near $(1,0)$. Note the accumulation of points near $(0.85,0)$,
$(0.7,0.6)$, $(0.4,0.4)$ and $(0.35,0)$, which are near the set $\Lambda$ and
its images under $\Phi$, as described above. Three initial points were used in
order to get a reasonable representation of the chaotic iterates because of
the sensitivity (associated with chaos) of the system and the limits of
computing accuracy.

\begin{figure}[th]
\centering
\includegraphics[width=5.5in]{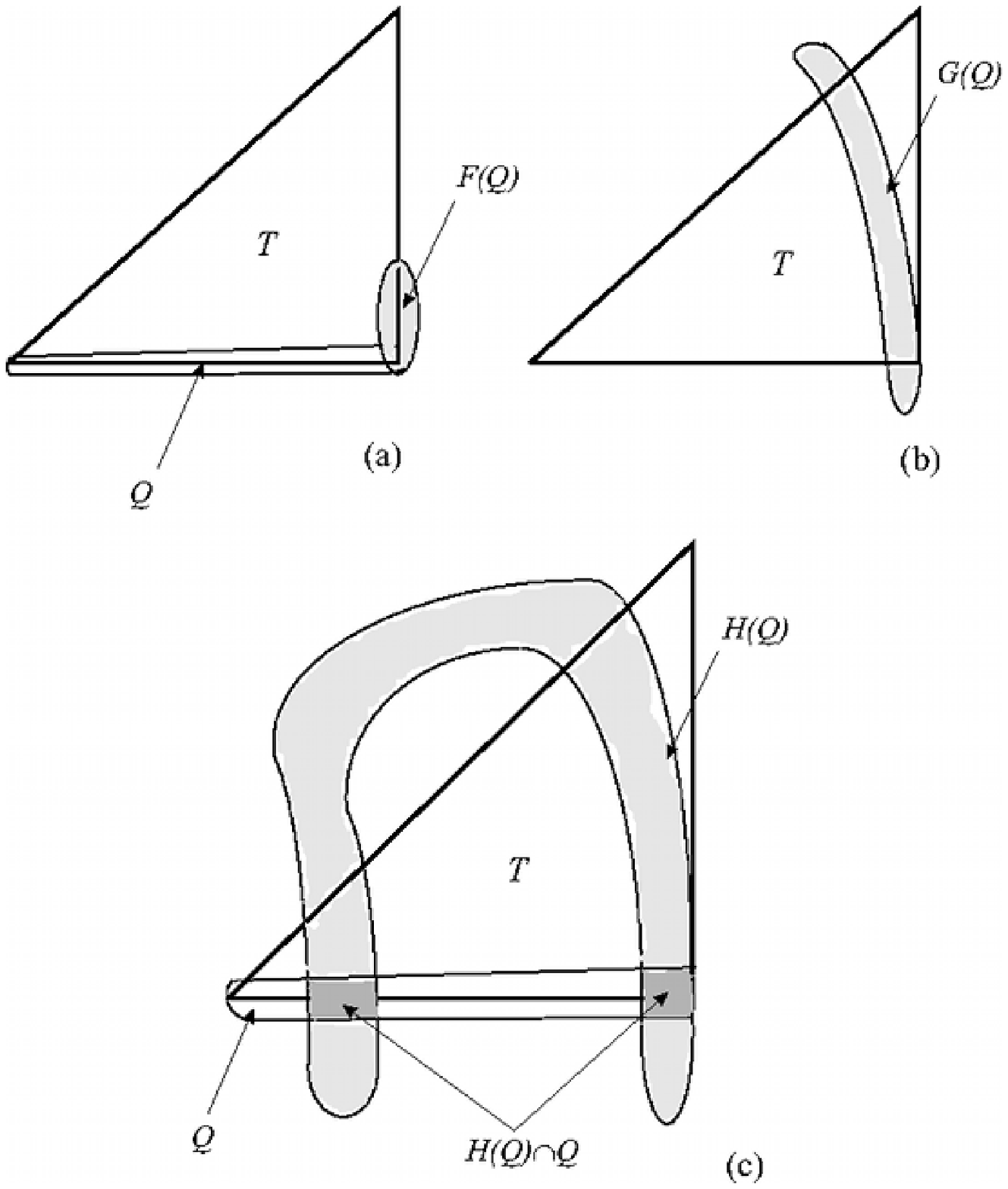} \caption{Horseshoe-like
behavior of map for $\lambda=0.99,\, \mu=5$: (a) Domain and first stage of
geometric image of iterated map; (b) Second stage focusing on stretching and
nascent folding; and (c) Final stage of horseshoe configuration}%
\label{ffchaos}%
\end{figure}

\begin{figure}[th]
\centering
\includegraphics[width=4.5in]{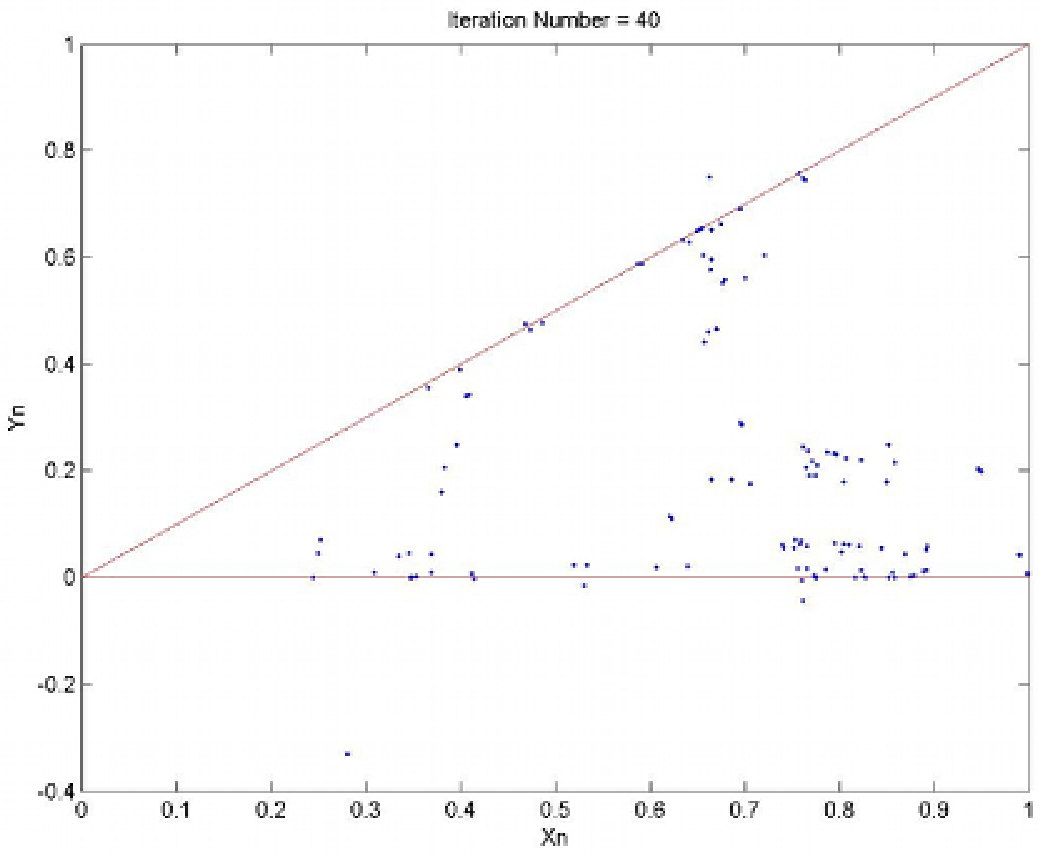}\caption{Chaotic orbits about
interior fixed point ($\cong(0.5569, 0.3569)$) for $\lambda=0.99,\, \mu=4.55$,
with eigenvalues $\xi\cong-0.5144\pm0.9596i$. The three initial points were
taken very close to $(1,0)$}%
\label{chaos}%
\end{figure}

\section{Comparison with Physical Models}

Our purpose in this section is to show that our discrete dynamical model
shares many properties with actual physical realizations (and their associated
mathematical models) of the \emph{R}-\emph{S} flip-flop circuit. Among the
several physically based studies of flip-flop type circuit behavior, which
includes the work in \cite{Chaney, Danca, KA, kac, LMH, Moser, msd, okt, rzw,
ztkbn}, perhaps the best source of comparison is provided by the investigation
of Okazaki \emph{et} \emph{al}. \cite{okt}, so this shall be our focus here.
We shall also provide additional numerical simulation illustrations of the
behavior of the orbits of our system to further highlight the areas of
agreement between our dynamical model and that of \cite{okt}.

The approach in \cite{okt} is quite different from ours. Okazaki and his
collaborators start with a physical realization of an \emph{R}-\emph{S}
flip-flop circuit using two capacitors, one inductor, one linear resistor, one
DC battery and a pair of piecewise linear resistors comprised of tunnel
diodes. They then derive the state equations for their circuit (hereafter
referred to as the OKT system), which is a system of three first-order,
piecewise linear ordinary differential equations in two dependent voltage and
one dependent current variable, and depends on several parameters associated
with the various electrical elements in the circuit. In their study they find
numerical solutions of this system of equations using a Runge-Kutta scheme,
which they compare with data collected directly from the physical circuit
using fairly standard electronic measuring and representation devices. They
find that the dynamical behavior deduced from numerical simulation of their
system of differential equations is in very good agreement with that observed
experimentally from the physical circuit. This behavior includes Hopf
bifurcation and chaotic dynamics for certain ranges of their parameters ,
which of course we have also observed and actually proved for our discrete
dynamical model.

More specifically, they deduce - mainly by observing the phase space behavior
of their numerical solutions - that there is a certain parameter value where
their system develops a Hopf bifurcation on a center manifold, followed, as
the parameter increases, by what appears to be a sequence of Hopf bifurcations
on the periodic orbits generated. This, of course, is consistent with the
qualitative behavior for our model described in Theorems 4.1 and 5.1, and
illustrated in Figs. 3 and 4. Moreover, although there are one-dimensional but
no two-dimensional Poincar\'{e} sections studied in \cite{okt}, several
projections of phase space structure on the two-dimensional voltage coordinate
plane, which are at least indicative of Poincar\'{e} map behavior on this
plane, have a striking similarity to the ring structure delineated in Theorem
5.1. Thus, there appears to be considerable qualitative agreement in the
dynamical behavior of our model and the OKT system for parameter values nearly
up to but just below those producing full-blown chaos.

Chaotic dynamics for the OKT system \cite{okt} is inferred primarily by
numerical computation of Lyapunov exponents, analysis of approximate
one-dimensional Poincar\'{e} sections, and observation of very complicated,
ostensibly random outputs in the experimental monitors. The projections onto
the voltage coordinate planes also have a somewhat chaotic appearance, with a
complicated looking tangle of orbits attached to and partially surrounding the
ring configuration mentioned above. Comparing this with the dynamics in the
examples pictured in Fig. 5 and 6, which shows ring configurations and the
tell-tale splatter (around the rings and concentrated near the fixed point)
associated with chaotic discrete dynamics, we have additional qualitative
validation for our model.

\section{Concluding Remarks}

In this investigation we have introduced and analyzed a rather simple discrete
dynamical model for the \textit{R}-\textit{S} flip-flop circuit, which is
based upon the iterates of a two-parameter family of planar, quadratic maps.
We proved that for certain parameter ranges, the dynamics of our model
exhibits the qualitative behavior expected in and observed for physical
realizations of the logical \textit{R}-\textit{S} flip-flop circuit such as
Hopf bifurcations and chaotic responses including oscillatory outputs of
arbitrarily large periods concentrated around states corresponding to nearly
equal set and reset inputs. In addition, we indicated how the dynamics of our
model displays fascinating complexity - as the parameters are varied -
generated by cascading bifurcations of stability transferring, doublings of
invariant collections of curves encircling a single equilibrium state, which
produce extremely intricate orbit structures.

Not being satisfied with the fact that the interesting variety of complex
dynamical structures produced by our model is certainly interesting from a
purely mathematical perspective, we undertook a more comprehensive validation
by comparing our results with the numerically simulated and experimentally
observed characteristics of a fairly standard realization of an \textit{R}%
-\textit{S} flip-flop circuit comprised of linear elements such as inductors
and capacitors and piecewise linear components consisting of tunnel diodes. We
found that the qualitative agreement between the dynamics of our model and
that of physical realization is surprisingly good. Notwithstanding this very
favorable comparison, we are aware that it may be largely fortuitous. After
all, our model is formulated in an essentially \emph{ad hoc} manner that
relies heavily on intuition and a desire to obtain the simplest maps producing
the kinds of dynamics known to be exhibited in working \textit{R}-\textit{S}
flip-flop circuits. Consequently, we are in the near future going to revisit
this logical circuit and investigate others of its kind using a much more
direct, first principles oriented approach along the lines of the work of
Danca \cite{Danca} and Hamill \emph{et al}. \cite{hdj}. Of course, it would be
particularly satisfying if we are able to show that such an approach produces
essentially the same discrete model for the \textit{R}-\textit{S} flip-flop
circuit investigated here, which is something that we expect but naturally
remains to be seen. We also intend to formulate and prove a generalized
version of Theorem 5.1, which we expect to have numerous applications in our
envisaged program of developing discrete dynamical models for a host of
logical circuits.


\begin{thebibliography}{99}                                                                                               %


\bibitem {VIA}V. Arnold, Small denominators, I: mappings of the circumference
into itself, \emph{AMS Transl. Ser. 2} \textbf{46} (1965), 213-284.

\bibitem {db}D. Blackmore, New models for chaotic dynamics, \emph{Regular \&
Chaotic Dynamics} \textbf{10} (2005), 307-321.

\bibitem {CB}J. Champanerkar and D. Blackmore, Pitchfork bifurcations of
invariant manifolds, \emph{Topology and Its Applications} \textbf{154 }(2007), 1650-1663.

\bibitem {Chaney}T. Chaney, A note on synchronizer and interlock maloperation,
\emph{IEEE Trans. Computers} \textbf{C-28 }(1979), 802-804.

\bibitem {chua}L. Chua, Chua's circuit: Ten years later, \emph{IEICE Trans.
Fundamentals} \textbf{E77-A} (1994), 1811-1821.

\bibitem {CWHZ}L. Chua, C.-W. Wu, A. Huang and G.-Q. Zhong, A universal
circuit for studying and generating chaos - Part II: Strange attractors,
\emph{IEEE Trans. Circuits \& Systems} \textbf{40} (1993), 745-761.

\bibitem {DMP}M. D'Amico, J. Moiola and E. Paolini, Hopf bifurcation in
discrete-time systems via a frequency domain approach, \emph{Proc. IEEE COC
Conf}., St. Petersburg, Russia, 2000, pp. 290-293.

\bibitem {Danca}M-F. Danca, Numerical approximation of a class of switch
dynamical systems, \emph{Chaos, Solitons \& Fractals} \textbf{38 }(2008), 184-191.

\bibitem {guho}J. Guckenheimer and P. Holmes, \emph{Nonlinear Oscillations,
Dynamical Systems, and Bifurcations of Vector Fields}, Springer-Verlag, New
York, 1983.

\bibitem {hdj}D. Hamill, J. Deane and D. Jeffries, Modeling of chaotic DC/DC
converters by iterated nonlinear maps, \emph{IEEE Trans. Power Electron}.
\textbf{7} (1992), 25-36.

\bibitem {hamill}D. Hamill, Learning about chaotic circuits with SPICE,
\emph{IEEE Trans. Education} \textbf{36} (1993), 28-35.

\bibitem {hartman}P. Hartman, \emph{Ordinary Differential Equations}, $2^{nd}$
\emph{ed}., Birkh\"{a}user, New York, 1982.

\bibitem {IJ}G. Ioos and D. Joseph, \emph{Elementary Stability and Bifurcation
Theory}, Springer-Verlag, New York, 1981.

\bibitem {KA}T. Kacprzak and A. Albicki, Analysis of metastable operation in
\textit{RS} CMOS flip-flop, \emph{IEEE J. Solid-State Circuits} \textbf{SC-22}
(1987), 57-64.

\bibitem {kac}T. Kacprzak, Analysis of oscillatory metastable operation of
\textit{RS} flip-flop, \emph{IEEE J. Solid-State Circuits} \textbf{23} (1988), 260-266.

\bibitem {kathas}A. Katok and B. Hasselblatt, \emph{Introduction to the Modern
Theory of Dynamical Systems}, Cambridge University Press, Cambridge, 1995.

\bibitem {LMH}G. Lacroix, P. Marchegay and N. Al Hossri, Prediction of
flip-flop behavior in metastable state, \emph{Electron. Lett}. \textbf{16}
(1980), 725-726.

\bibitem {MM}J. Marsden and M. McCracken, \emph{The Hopf Bifurcation and Its
Applications}, Springer-Verlag, New York, 1976.

\bibitem {Moser}J. Moser, Bistable systems of differential equations with
applications to tunnel diode circuits, \emph{IBM J. Res. Dev}. \textbf{5
}(1961), 226-240.

\bibitem {msd}K. Murali, S. Sinha and W. Ditto, Implementation of a nor gate
by a chaotic Chua's circuit, \emph{Int. J. Bifurcation and Chaos} \textbf{13}
(2003), 2669-2672.

\bibitem {okt}H. Okazaki, H. Nakano and T. Kawase, Chaotic and bifurcation
behavior in an autonomous flip-flop circuit used by piecewise linear tunnel
diodes, \emph{IEEE Proc. ?}, 1998, III-291 - III-297.

\bibitem {palmel}J. Palis and W. de Melo, \emph{Geometric Theory of Dynamical
Systems}, Springer-Verlag, Berlin, 1982.

\bibitem {rzw}M. Ruzbehani, L. Zhou and M. Wang, Bifurcation features of a
dc-dc converter under current-mode control, \emph{Chaos, Solitons \& Fractals}
\textbf{28 }(2006), 205-212.

\bibitem {wan}Y.H. Wan, Computations of the stability condition for the Hopf
bifurcation of diffeomorphisms on $\mathbb{R}^{2}$, \emph{SIAM J. Appl. Math}.
\textbf{34} (1978), 167-175.

\bibitem {WXH}G. Wan, D. Xu and X. Han, On creation of Hopf bifurcations in
discrete-time nonlinear systems, \emph{CHAOS} \textbf{12} (2002), 350-355.

\bibitem {wigbook}S. Wiggins, \emph{Introduction to Applied Nonlinear
Dynamical Systems and Chaos}, \emph{2}$^{nd}$\emph{\ ed}., Springer-Verlag,
New York, 2003.

\bibitem {ztkbn}A. Zorin, E. Tolkacheva, M. Khabipov, F.-I. Buchholz and J.
Niemeyer, Dynamics of Josephson junctions and single-flux-quantum networks
with superconductor-normal-metal junction shunts, \emph{Phys. Rev. B}
\textbf{74} (2006), 014508.
\end{thebibliography}
\end{document}